\documentclass[a4paper]{article}
\usepackage[margin=0.9in]{geometry} %

\usepackage{amsmath,amssymb,amsfonts}
\usepackage{algorithmic}
\usepackage{graphicx}
\usepackage{textcomp}
\usepackage{xcolor}
\usepackage{hyperref}
\usepackage{xspace}
\usepackage{comment}
\usepackage{float}
\usepackage{todonotes}            %
\usepackage[symbol]{footmisc}

\def\orcidID#1{\unskip$^{[#1]}$} %

\begin{document}

\title{Run your HPC jobs in Eco-Mode: revealing the potential of user-assisted power capping in supercomputing systems}
\author{Luc Angelelli\orcidID{0009-0008-2792-6264} \and Danilo Carastan-Santos\orcidID{0000-0002-1878-8137} \and Pierre-François Dutot\orcidID{0000-0002-1992-3388}}
\date{Univ. Grenoble Alpes, CNRS, INRIA, Grenoble INP*, LIG \\38000
 Grenoble, France}

\maketitle

\begin{abstract}

The energy consumption of an exascale High-Performance Computing (HPC) supercomputer rivals that of tens of thousands of people in terms of electricity demand.
Given the substantial energy footprint of exascale HPC systems and the increasing strain on power grids due to climate-related events, electricity providers are starting to impose power caps during critical periods to their users.
In this context, it becomes crucial to implement strategies that manage the power consumption of supercomputers while simultaneously ensuring their uninterrupted operation.

This paper investigates the proposition that HPC users can willingly sacrifice some processing performance to contribute to a global energy-saving initiative.
With the objective of offering an efficient energy-saving strategy by involving users, we introduce a user-assisted supercomputer power-capping methodology.
In this approach, users have the option to voluntarily permit their applications to operate in a power-capped mode, denoted as 'Eco-Mode', as necessary.

Leveraging HPC simulations, along with energy traces and application metadata derived from a recent Top500 HPC supercomputer, we conducted an experimental campaign to quantify the effects of Eco-Mode on energy conservation and on user experience.
Specifically, our study aimed to demonstrate that, with a sufficient number of users choosing Eco-Mode, the supercomputer maintains good performances within the specified power cap.
Furthermore, we sought to determine the optimal conditions regarding the number of users embracing Eco-Mode and the magnitude of power capping required for applications (i.e., the intensity of Eco-Mode).

Our findings indicate that decreasing the speed of jobs can decrease significantly the number of jobs that must be killed.
Moreover, as the adoption of Eco-Mode increases among users, the likelihood of every job to be killed also decreases.

\vspace{2mm}
\end{abstract}

\section{Introduction}
\label{sec:introduction}

\footnotetext[1]{Institute of engineering Univ. Grenoble Alpes}

High Performance Computing is being essential for advancements in many fields of science and engineering.
However, operating current exascale machines requires substantial amounts of power.
For instance, the current first-ranked supercomputer in the Top500 list~\cite{top500} consumes 29 megawatts~\cite{frontierspecs} of power at peak performance.
Such a power consumption is equivalent to the electricity demand of around twenty thousand US citizens, taking into account an average 
US \textit{per capita} electricity demand of 1285W.\footnote[2]{Source: \url{https://en.wikipedia.org/wiki/List_of_countries_by_electricity_consumption}}

This energy consumption puts a heavy strain in the electricity grid.
This strain increases the risk of power outages, especially in the advent of increased electricity demand induced by climate events.
An example is the 2021 Texas power crisis~\cite{texaspoweroutage}: severe winter storms caused a state-wide power outage due to increased heating demand.

To properly react in face of similar power crises, exascale supercomputer maintainers must employ measures to control the power consumption of their platforms.
A typical measure is to enforce a power cap~\cite{borghesi2015powercap,kontorinis2012powercap} in the platform when needed.
Power-constrained supercomputers introduces the challenge of how to deal with the applications to both comply with the power cap and maintain a sufficient operational level.
One challenge is therefore on how to deal with the running applications to comply with the power cap.

Typical approaches to deal with the power cap are (i) applying a Dynamic Voltage and Frequency Scaling (DVFS) to all running applications~\cite{nana2023energyopti} or (ii) killing running applications until the power cap is satisfied~\cite{maiterth2018jobkilling,nana2023energyopti}.

We explore the idea that users may volunteer and allow their applications to run
in a lower power state when needed. This lower power state will result in
increased processing time for the selected applications. However, volunteering
to run the applications in a lower power state may reduce the chances of quality
of service (QoS) disruptions for the user, specifically in the case where
running jobs are killed to comply with the power cap. The main research question
tackled in this paper is therefore: \textit{facing enforced power reductions in
a supercomputer, should users opt to run their applications in a lower power
state to reduce the chances of their applications to be killed?}

 We propose an experimental method that exploits simulation and power
 consumption traces of real-world supercomputers to predict the outcomes of
 several power cap situations and user's lower power acceptance levels. With our
 method, we performed a detailed experimental campaign to understand the
 aforementioned trade-offs. In short, we show experimental evidence that
 highlight the following findings.

\begin{itemize}
  \item If 100\% of the users accept to run applications in a lower power state,
  we can comply with moderate levels of power cap (up to 70\% of the maximum
  power in our experiments) with no job kills. As the level of acceptance
  reduces, the number of job kills increases to comply with the power cap, thus
  degrading the QoS of the platform.
  \item Volunteering to run applications in a lower power state does not affect
  the platform performance in handling the workload.
\end{itemize}

The remaining parts of this paper are organized as follows. We present an
overview of the related works in Section~\ref{sec:methods}, and preliminary
definitions in Section~\ref{sec:preliminary-concepts}. We present our proposed
methodology in Section~\ref{sec:methods}, present the experimental results in
Section~\ref{sec:results}, and we finally conclude and discuss future works in
Section~\ref{sec:conclusion}.

\section{Related Work}
\label{sec:related-work}

Below, we briefly present recent research trends and related works relevant to power/energy-aware HPC resource monitoring and management.
We invite the reader to consult Maiterth \textit{et al.}~\cite{maiterth2018energy} and Kocot \textit{et al.}~\cite{kocot2023energy} for a more detailed survey on the subject.

Researchers have been proposing energy/power aware scheduling methods for HPC platforms by employing a large variety of methods, ranging from integer programming~\cite{datazero} to heuristics~\cite{chasapis2019power,hu2021characterization}, and Machine Learning~\cite{d2021energy}.
Most of these works calculate the power consumption of an application relying on the Thermal Design Power (TDP) of the processors/accelerators, and a mix of the application's processing time and resource utilization.
Several works also use energy measurements from interfaces such as RAPL~\cite{khan2018rapl,saurav2016rapl}, but in aggregations such as average energy/power consumption.
In contrast, our work relies on real-measured data consisting of fine-grained time series of the power consumption of the applications to show the efficiency of our method.
Few works rely on time-series data of the power consumption~\cite{chasapis2019power,d2021energy}, and it is mainly exploited in studies that analyze the behavior of the platforms~\cite{patel2020does,shin2021revealing}.
This work is a step towards exploiting power consumption time-series data for power-aware HPC resource management.

In summary, the majority of works present sophisticated methods for energy-aware resource management.
However, it is widely known that sophisticated methods hinders their deployment in real-world production platforms~\cite{feitelson1997theory}.
Furthermore, these sophisticated methods may require additional \textit{a priori} characterizations~\cite{chiesi2014power} or predictions~\cite{borghesi2016predictive,chasapis2019power,d2021energy,frey2022benchmarking} of the power consumption of the applications.
Predicting the power consumption can be quite hard to achieve.
For instance, in the case of GPUs, some applications can result in high power consumption variability when running in distinct GPUs of the same model and vendor~\cite{sinha2022not}.

Our work is in the category of power capping for power-aware resource management in power constrained platforms~\cite{etinski2012parallel,georgiou2015adaptive,borghesi2018scheduling,zhao2023sustainable}.
Power capping consists in adjusting the computing nodes to consume less power than its allowed maximum by employing, for instance, Dynamic Voltage and Frequency Scaling (DVFS).
We distinguish ourselves from previous works in power capping because we rely on the users' engagement to run applications under power cap.
Engaging users results in a simple power-aware resource management method that does not need any \textit{a priori} characterization or prediction of the power consumption of the applications, making it easier to deploy in production platforms.
The principle is straightforward: we can comply with the power cap if enough users volunteer to run their applications under slight levels of power cap.
No need for complex power-capping decision-making from the side of the platform.
We contribute with an experimental campaign to show evidence of the efficiency of such a simple method.

\section{Preliminary definitions}
\label{sec:preliminary-concepts}

We consider an HPC platform as a computing cluster with $m$ computing nodes that are connected by a certain interconnection topology.
Each computing node has one or more CPUs and also accelerators (e.g., GPUs).
In the context of this paper, we consider the accelerators as only GPUs.
The computing clusters are homogeneous in the context of the computing nodes.
That is, all computing nodes have the same number and model of CPUs and GPUs.

During the platform's operation, several users submit applications, hereafter referred to as jobs.
To deal with the jobs' submission and processing in the platform, a management system called Resources and Jobs Management System (RJMS) runs in the platform.
The RJMS is the main interaction point between the users and the HPC platform.
Users can submit applications to the RJMS at any point in time, and the RJMS has no ahead information about which jobs will arrive.
The literature refers to this job submission characteristic as \textit{online} job submission.
We assume that the RJMS assigns the processing order of the jobs in First Come First Served (FCFS) order.
FCFS is the baseline job queue ordering heuristic in the popular backfilling scheduling algorithm~\cite{mu2001utilization}.

An electricity provider powers the platform with electricity. We consider that
there are certain time periods where the electricity provider is not able to
fully power the platform, resulting therefore in power capped periods during the
platform's operation. In the real-world, this cap can come from sources such as
(i) general electricity grid overdemand due to an external event, and (ii)
temporary reduction of the renewable power capacity of the electricity provider.
The electricity provider informs the platform maintainers in advance about
when and how much the HPC platform will be power capped to accommodate the
peaks of demand.

The RJMS is also responsible to manage the jobs that are currently running in
the HPC platform. Therefore, once a power capped period arrives, it is up to the
RJMS to deploy measures to comply with the power cap. We elaborate in these
measures in the section below.

\section{Proposed methodology}
\label{sec:methods}

\subsection{Methods to comply with a power cap}

In production clusters such as Marconi100, the scheduling algorithm used is
usually a variant of FCFS, or EASY-bf~\cite{mu2001utilization}. As such, in the
need to set a power cap, production clusters can adopt a variant of their
preferred scheduling algorithm that would either kill jobs arbitrarily, or make
a reservation for machines to shut down, such that the maximum power draw of the
remaining machines is below the powercap. We evaluate in this paper two
scheduling algorithms: 

\textit{FCFS\_killer} (Baseline): kills jobs until the power cap is reached,
newest job first. Upon job termination, the RJMS can shut down the computing
nodes previously allocated by the killed job, thus saving power. This method is
arguably the easiest to be deployed in practice, but it can be very intrusive,
affecting the Quality of Service (QoS) of the platform, since
\texttt{FCFS\_killer} risks loosing all the computation that was being performed
by the killed job. We decided to select the newest job first in an attempt to
reduce the wasted energy caused by killing.

\textit{FCFS\_eco\_mode} (Our contribution): on the job submission, users may
flag their jobs, indicating that the job can be run in a slowed down state. We
hereafter refer to this flagged job as EcoJob. During the power-cap period, the
RJMS looks for flagged jobs that are running, and uses DVFS to proportionally
slow down the nodes assigned to these flagged jobs, until the power cap is
reached. The RJMS slows down the assigned nodes up to a lower limit of 50\% of
the maximum power of the node. If all possible nodes are slowed down to this
limit and the power cap is still not satisfied, a \texttt{FCFS\_killer} routine
starts, which will kill jobs until the power cap is satisfied. Non-flagged jobs
will be prioritized to be killed by the \texttt{FCFS\_killer} routine.

There are 4 events that the schedulers need to react to in a power-capping setting.
\texttt{FCFS\_eco\_mode} reacts to these events as follows:

\begin{algorithmic}
    \REQUIRE{{Powercap Period is Incoming}}
        \FORALL{EcoJob currently running}
            \STATE set DVFS to the lowest mode
        \ENDFOR
        \IF{current power $>$ powercap}
            \STATE kill jobs, non-EcoJobs first, newest first
        \ELSIF{current power $<$ powercap}
            \STATE increase DVFS until powercap is met
        \ENDIF

    \REQUIRE{{Powercap Period Ended}}
        \FORALL{EcoJob currently running}
            \STATE set DVFS to the highest mode
        \ENDFOR

    \REQUIRE{{Job finished}}
        \STATE Liberate resources
        \STATE Reset corresponding machines' DVFS state
        \STATE Execute jobs from queue as long as they fit in power and resources

    \REQUIRE{{Job submitted}}
        \STATE Add to queue
        \STATE Execute jobs from queue as long as they fit in power and resources
        
\end{algorithmic}

The idea behind \texttt{FCFS\_eco\_mode} is to let the user decide to trade-off
performance (i.e., let their jobs run slower) to increase the chances that their
jobs will not be killed. \texttt{FCFS\_eco\_mode} tries to engage users as a
front line measure to comply with the power cap, potentially avoiding degrading
the QoS with job kills. It is important to emphasize that this idea behind
\texttt{FCFS\_eco\_mode} is agnostic for any kind of scheduling algorithm, since
it concerns the jobs that are already running.

Understanding the effects of \texttt{FCFS\_killer} and \texttt{FCFS\_eco\_mode} by experimenting them in real-world platforms is risky, onerous, and time-consuming.
We address this hindrance by performing simulation experiments of operating an HPC platform.
In the next section, we present the details to achieve in this simulation.

\subsection{Simulating a supercomputer under power cap}

\subsubsection{Data sources for simulation}

It is a common practice from supercomputer maintainers to register information about the operation of the platform, notably details about the jobs in the form of workload traces (e.g., number of requested processors, arrival time, run time).
An example is the workload traces present in the parallel workloads archive~\cite{feitelson2014experience}.
However, such traces do not contain information about the power consumption of the computing nodes during the platform's operation, which hinders the task of simulating the workloads traces taking into account the energy consumption and power cap.

To overcome this hindrance, we exploit a recent dataset from the Marconi100 supercomputer~\cite{borghesi2023m100}.
This dataset not only contains the aforementioned information about the jobs, but it contains as well time-series data about the power consumed by the computing nodes.

On top of the common data collection for job information in a Job Table, two plugins were used to collect energy information, IPMI and Ganglia.
The IPMI (Intelligent Platform Management Interface) plugin serves as a data collection tool, retrieving information from the Out-of-Band (OOB) management interface, specifically the Baseboard Management Controller (BMC) of the computing nodes.
The BMC is a hardware component embedded in servers or cluster nodes that facilitates remote monitoring and management of the system.
The IPMI plugin gathers sensor data -- notably the power consumption -- from the BMC installed in the nodes.
The Ganglia plugin serves as an energy monitoring tool for most components of each node.
In particular, it recorded the energy use of each GPU for every node, with a similar frequency to IPMI on the CPUs.

To establish at which points in time the HPC platform will be power capped, we
considered a use-case from the available electricity demand data from
RTE~\cite{RTE} (\textit{Réseau de Transport d'Electricité}). RTE's data gives us
real-time data and history on the electricity demands in France. With this data
we can observe regular daily peaks around 6PM to 8PM. Therefore, in our use-case
we suppose that the power cap periods happen every day between 6PM and 8PM.

\subsubsection{Preparing data for simulations}

\begin{figure*}[h!]  
  \centering \includegraphics[width=\linewidth]{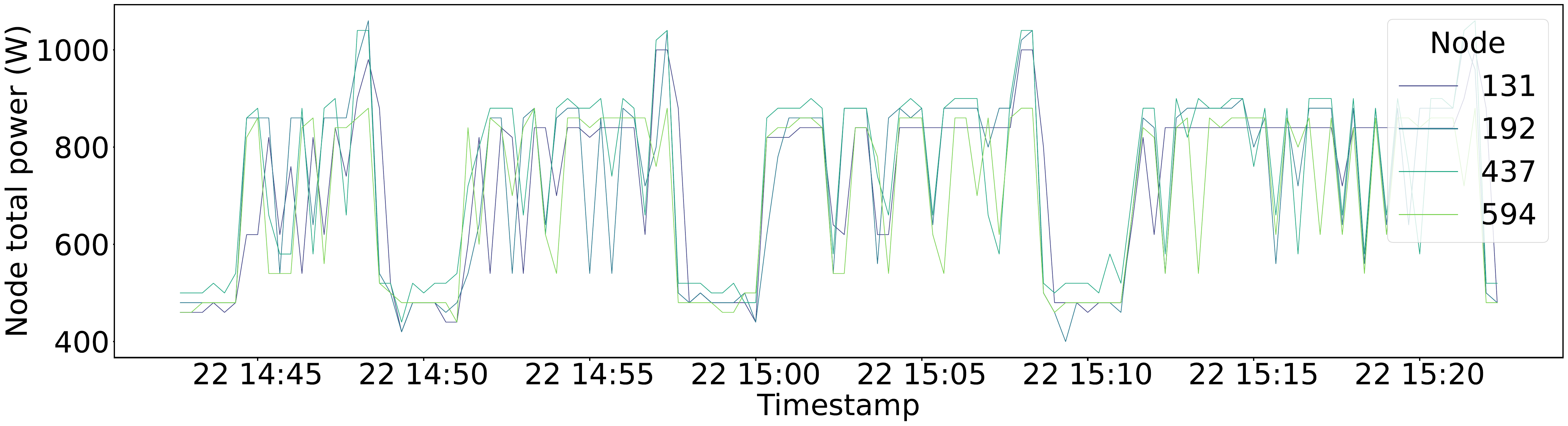}\\
  \centering \includegraphics[width=\linewidth]{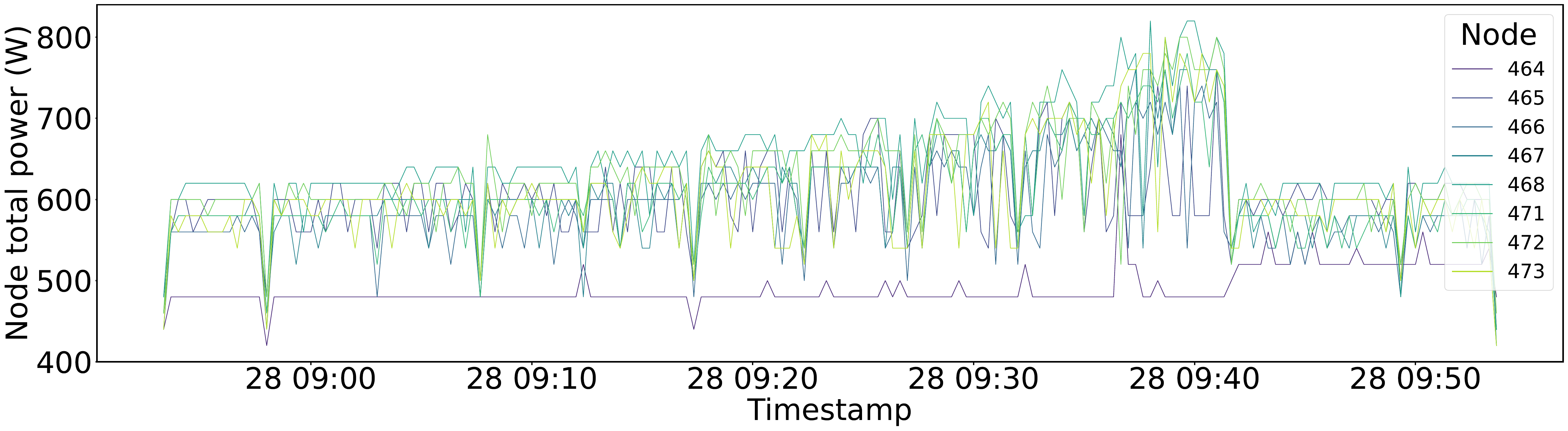}  
  \caption{Examples illustrating the power consumption profiles of nodes of two jobs, one job per graph, obtained from the Marconi100 dataset.
Each job is a multi-node application, i.e., four nodes in the top graph and eight nodes in the bottom graph.}
  \label{fig:m100-energy-profiles}
\end{figure*}

By crossing the timestamp information about the start/end time of the jobs (Job Table data) and the time-series data of the computing nodes' power consumption (IPMI Data and Ganglia Data), we can extract time-series data of the power consumption of the jobs.
Figure~\ref{fig:m100-energy-profiles} illustrates a few examples of the extracted jobs' power consumption data.
From this per-job power consumption data, we extract the power consumed by the CPUs and GPUs during the jobs' execution as a time-series data, sampled once every 20 seconds.

To simulate the jobs taking into account the sampled power data, we converted
the power profile into a computing profile of the jobs as follows: for a job $j$
and its power profile time-series data $P^{cpu}_j$ and $P^{gpu}_j$, we use
Equations~\ref{eq:conversion-cpu} and~\ref{eq:conversion-gpu} to calculate the
computing profile time-series data $C^{cpu}_j$ and $C^{gpu}_j$.

\begin{equation}
  C^{cpu}_j = \left( max(-602 + 134*log(P^{cpu}_j), 0) \right)
  \label{eq:conversion-cpu}
\end{equation}

\begin{equation}
  C^{gpu}_j = \left( max(-27.1 + 7.3*log(P^{gpu}_j), 0) \right)
  \label{eq:conversion-gpu}
\end{equation}

Those equations and the constant values used are derived from the CPU and GPU modeling for the Marconi cluster, as explained in the platform modeling paragraph below.
As such they allow us to simulate DVFS in a non-trivial fashion. The jobs are more than the simplistic rectangular job of constant resource use, and the machines have a more realistic behaviour thanks to the non-linearity of their model.

After this conversion, the jobs are now modeled as a sequence of amounts of
computing to be processed for every 20 seconds. This new modeling of the jobs
constitutes a more fine-grained computing profile of jobs, which gives us
flexibility to simulate the DVFS impact in the jobs in precise time windows
during the jobs' execution, and therefore to simulate the effects of a power cap on the jobs.
The power state employed for task execution is contingent upon its
inherent characteristics. Tasks characterized by significant
input/output (I/O) operations exhibit differing power state
requirements compared to those involving intensive numerical
computations. Consequently, the execution time of I/O-bound tasks
experiences less pronounced sensitivity to reduced clock frequencies,
in contrast to computationally intensive tasks.

\subsubsection{Configuring the simulations}

We use Batsim~\cite{Batsim} for all our simulations.
As a simulator, Batsim is able to modulate parameters such as frequency and power during simulation, mimicking a DVFS behavior.
For Batsim to work, it needs a platform, a workload, and a scheduler.

The platform has been made to resemble Marconi100.  As such, it has 5
batsim computing nodes, 1 representing the CPUs and 4 representing the
GPUs, for each real Marconi100 nodes.
The lack of publicly available studies on Marconi100 components and
the inherent limitations of production environments in revealing job
details necessitate assumptions for converting Marconi100 traces to a
batsim workload. As we lack most of the specifics of each job, we
decided to ignore communications, and focus only on computations.
Fortunately, analysis in~\cite{zacharov_zhores_2019}
by Zacharov \textit{et al.} enables direct conversion of power consumption to
GFlop/s for GPUs. However, no such straightforward method exists for
CPUs. Therefore, we scaled a CPU with known DVFS power consumption to
match the idle and maximum power of Marconi100 CPUs.
This approach yields a batsim platform emulating Marconi100, featuring
an approximation of DVFS for both CPUs and GPUs.

In our scenario, we imagine that the electricity provider imposes a power cap on the platform at specific times.
This power cap is given as dates in a configuration file, specifying for each change the time and the new power cap.
This power cap information is then used by the scheduler to place, slow or kill jobs when needed.

We exploit Batsim's composed jobs representation to use the calculated computing
profile time-series data (Equations~\ref{eq:conversion-cpu}
and~\ref{eq:conversion-gpu}) as workload inputs for the job submission into the
simulation.
To fairly compare between all the schedulers and parameters, we
dynamically submitted jobs to the platform to simulate users behavior
for a set amount of time (10 days in our experiments).  Submitting
jobs dynamically instead of replaying a static workload with fixed
submission times is important, as it allows us to regulate the job
queue, preventing some algorithms from having the unfair advantage of
being able to select from a large pool of jobs simply because they
mismanaged the beginning of the schedule and created a large backlog.

Additionally, comparisons between different power cap settings would
yield limited insights if the workload remained static across all
scenarios. Replaying the same job submissions under varying power caps
would primarily reflect the available energy, rather than provide
meaningful differentiation based on the evaluated parameters.

\section{Experimental Results}
\label{sec:results}

\subsection{Experimental Setup}

In our experiments, we compare two schedulers, \texttt{FCFS\_killer} and \texttt{FCFS\_eco\_mode} (Section~\ref{sec:methods}).
\texttt{FCFS\_killer} acts as our control.
It is a modified version of a standard FCFS algorithm, that reacts to a power cap event by killing jobs, most recent first, until the powercap is reached.
\texttt{FCFS\_eco\_mode} is also based on FCFS, but instead of killing jobs as soon as a power cap arrives, it tries to slow the EcoJobs (i.e., jobs that are flagged to run in a lower power state) first.
If slowing is not enough, then it kills jobs. Non-EcoJobs and newest jobs are selected first to be killed.
So killing jobs is its last resort.
With this approach, we hope to reduce the energy lost in killing jobs, reducing waste and improving quality of service.

There are 3 main parameters for our experiments.
Those are (i) the percentage of EcoJobs in the workload, (ii) the amount of power available when in power cap, and (iii) the scheduling algorithm used.

For the percentage of EcoJobs, we decided to use 0, 10, 25, 50, 75 and 100 percent to give us a good idea of the effects of adherence on the results.

The amount of power available when in power cap was set between 50 and 100 percent of the maximum power of the computing nodes.
We decided to not go below 50 percent because the nodes actually consume power even at idle, and reducing that power necessitate a different approach, like machine extinction~\cite{dutot2017extinction}.

\begin{figure}[h!]
    \centering
    \includegraphics[width=0.8\textwidth]{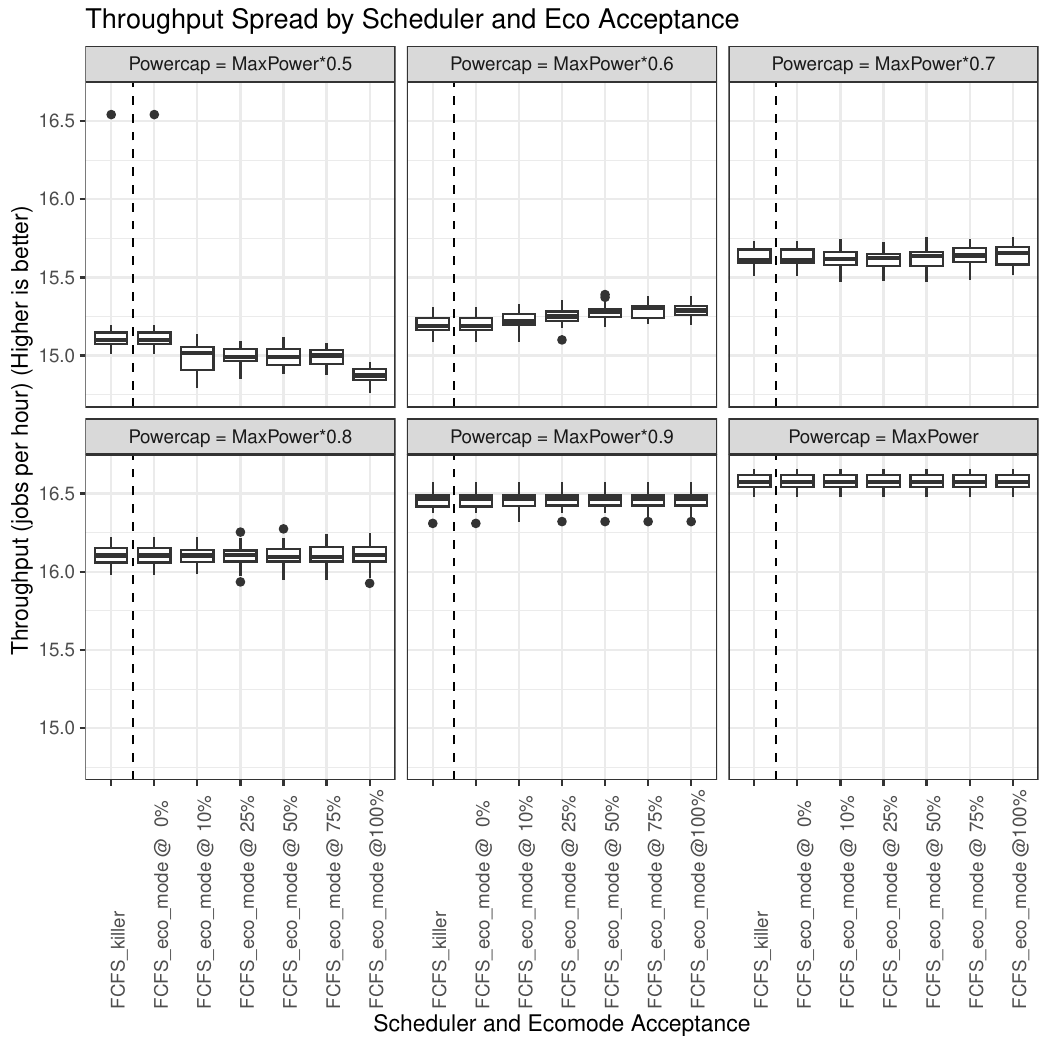}
    \caption{Throughput depending on the scheduler, the percentage of EcoJobs and the Strength of the power cap (i.e.
\% of maximum power available in a power cap period)}
    \label{fig:P_Experiments_Throughput-by-scheduler}
\end{figure}

\begin{figure}[h!]
    \centering
    \includegraphics[width=0.8\textwidth]{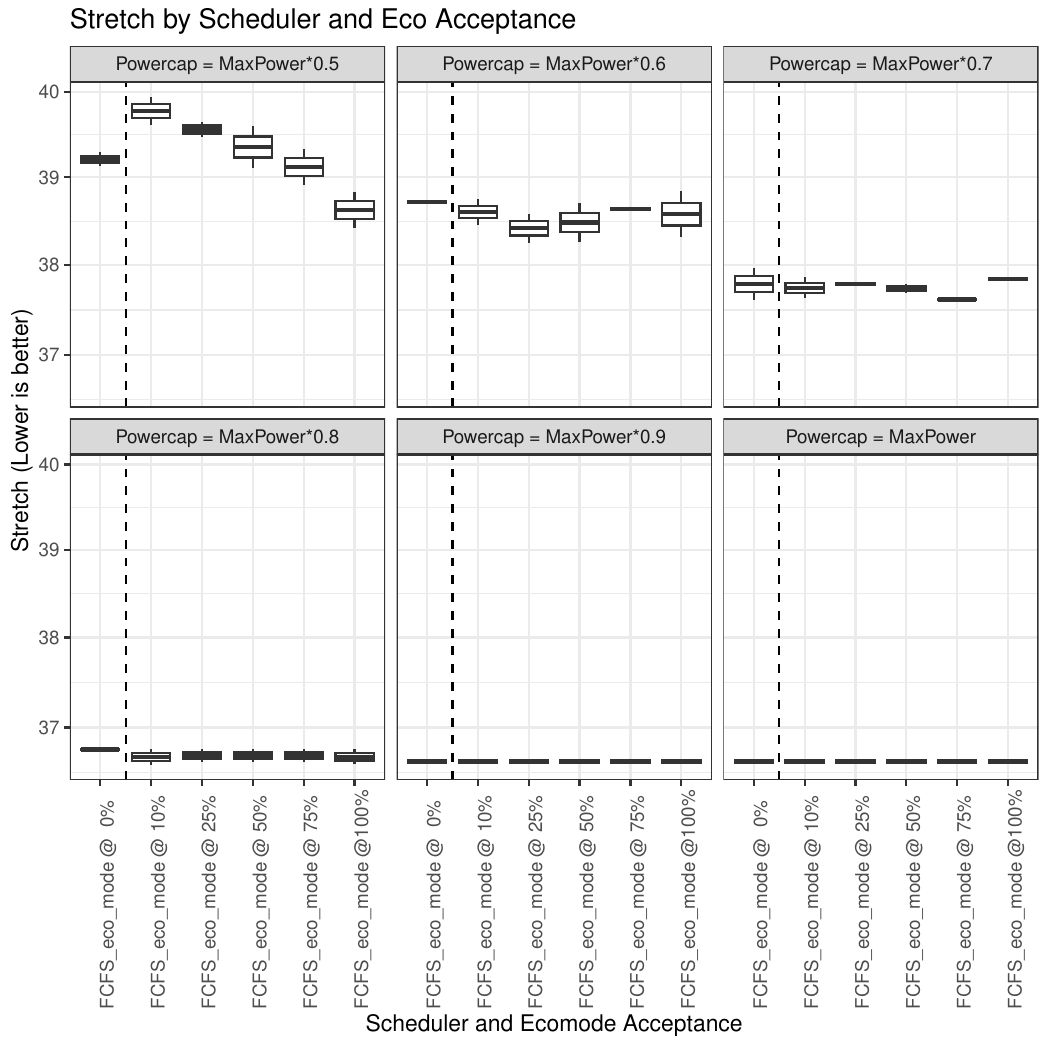}
    \caption{Stretch depending on the scheduler, the percentage of EcoJobs and the Strength of the power cap}
    \label{fig:P_Experiments_Stretch-without-kills}
\end{figure}

\begin{figure}[h!]
    \centering
    \includegraphics[width=0.8\textwidth]{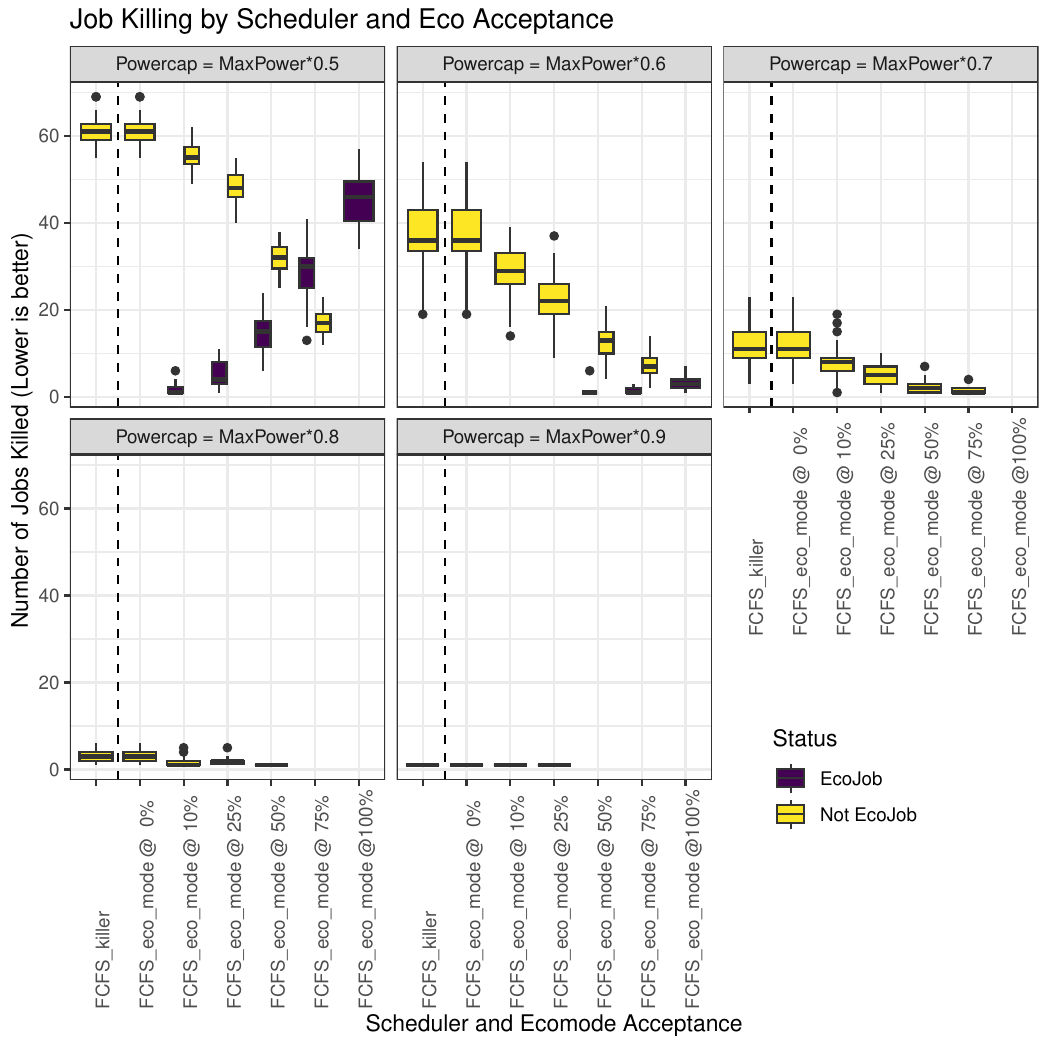}
    \caption{Number of jobs killed depending on the scheduler, the percentage of EcoJobs and the Strength of the power cap}
    \label{fig:P_Experiments_Kills-by-scheduler}
\end{figure}

\begin{figure}[h!]
    \centering
    \includegraphics[width=0.8\textwidth]{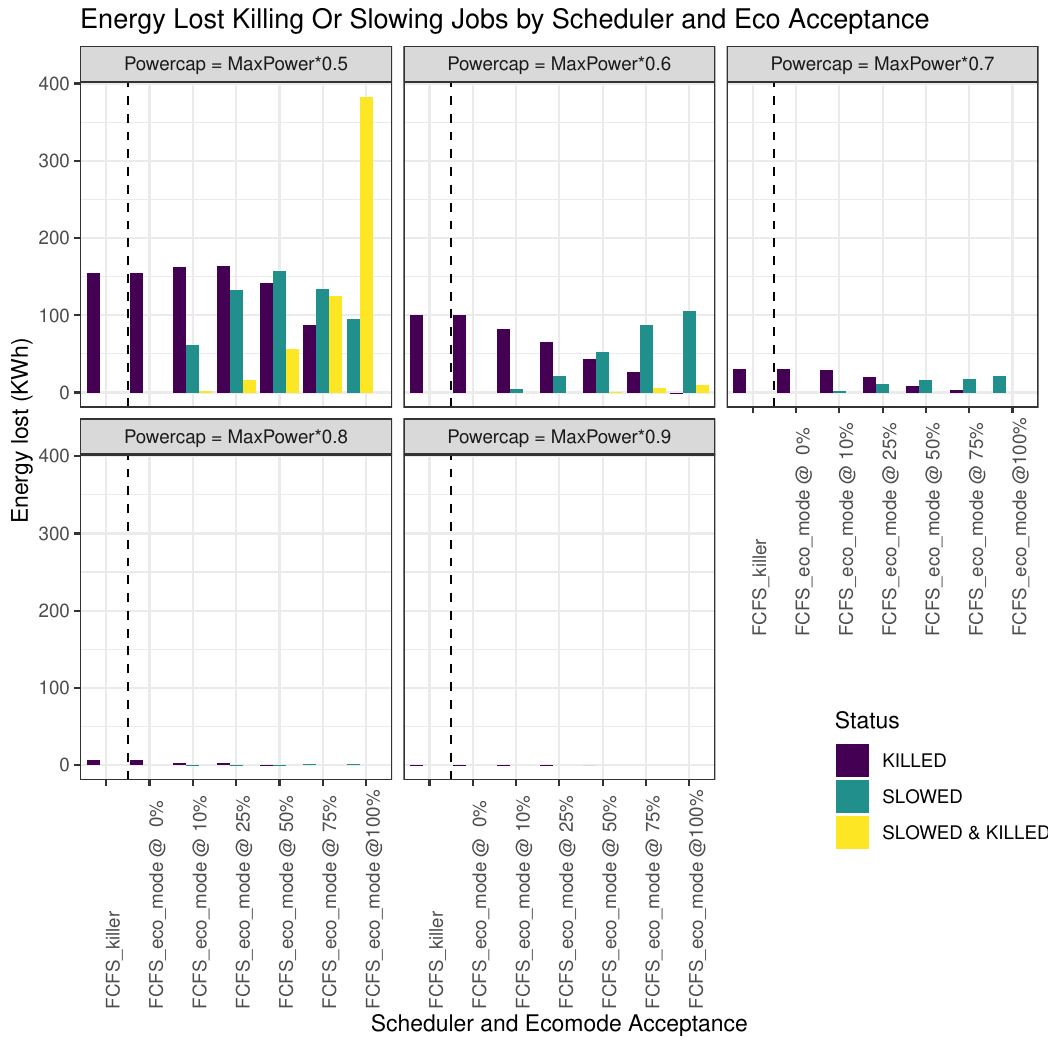}
    \caption{Energy lost by killing or slowing jobs depending on the scheduler, the percentage of EcoJobs and the Strength of the power cap}
    \label{fig:P_Experiments_Lost-energy-by-scheduler}
\end{figure}

\subsection{Results}

Our experiments consist of 30 independent runs for each combination of parameters.
Each run uses a seed equals to its run number.
The randomness in these experiments is used to define the order of arrival of tasks in the queue.
This allows us to see a wide range of behavior for the algorithms while still being able to compare equivalent runs between categories.

Each run is simulating 10 days of use of an 8 nodes machine, with each node being one CPU and 4 GPUs, with a power cap active between hours 18 and 20 each day. This platform is only a subset of the Marconi platform, to allow for a large parameter sweep. Larger experiments were done with more nodes, up to 98 (to represent 10\% of the Marconi platform). Those results were consistent with the results presented here, but we chose to present the full parameter exploration at a smaller scale rather than focus on a smaller set of results (as the running time of such experiments is quickly getting expensive in CPU time).
Let us first evaluate the performances of the algorithms.
We first observed the average throughput in job per hour, as in the number of jobs completed successfully within the 10 days, divided by the total number of hours.
We can see in Figure~\ref{fig:P_Experiments_Throughput-by-scheduler} that the average throughput is barely affected, less than 2\%, by our algorithm within each level of powercap.
The most prominent feature of these graphs is actually the restriction imposed by the power cap, accounting for about 10\% change for the extreme case of 50\% of power available only between 18 and 20. This difference is much larger than the $~4$\% difference in energy, since some jobs can be killed and jobs running at a lower power setting will generally consume more energy. 
We can also observe on this figure and other comparison of schedulers that \texttt{FCFS\_eco\_mode} with no adherents is exactly the same as \texttt{FCFS\_killer}.
This is no coincidence, as EcoMode was designed to be a sort of add-on to FCFS.

Observing the average stretch in Figure~\ref{fig:P_Experiments_Stretch-without-kills} yields a similar story, a very small effect from the algorithm, about 2\%, and a larger difference due to the power capping, about 5\%.
With these observations, we can say with some confidence that our algorithm is performing about the same as our control.
Finally, in Figure~\ref{fig:P_Experiments_Kills-by-scheduler}, we show the main arguments for choosing our algorithm over an FCFS variant, which is the amount of jobs killed to ensure powercap.
We show a consistent and significant decrease in the number of kills as the acceptance of ecomode increases.
Furthermore, our algorithm prioritize non-EcoJobs jobs when it needs to kill, so EcoJobs are way less likely to be sacrificed for the good of the platform.

The decrease in number of kills leads to a decrease in the energy wasted by killing, as shown in Figure~\ref{fig:P_Experiments_Lost-energy-by-scheduler}.
This energy is converted almost exactly into energy used by slowed jobs when the power cap is not too strong.
When the powercap is too strong, at 0.5 for example, we can see that tasks can be slowed, and then killed, leading to higher wasted energy than normal.

This is due, in part, to what we consider a waste when killing and when slowing jobs.
When killing jobs, we lose all the energy we put in calculation for that job.
But when we slow a job, because it takes longer to complete and the power to performance ratio of the compute nodes are non-linear, the job can use more energy to complete.
So we consider a waste the difference in energy between a slowed job and its full-speed counterpart.

If the job finishes, we could argue that it was not actually a waste of energy, as the alternative was to kill to make space for full-speed jobs.
But when we kill an already slowed job, we have to count it as a kill, and all energy spent on it is lost.

There are ways to mitigate these effects, like preemptively killing jobs, but we did not explore them in this paper.
But the more important point is that, if we do not turn off machines, the power cap can be too strong for any algorithm, as there is a power cap beyond which no job can be scheduled.
There is also a minimum power cap where all machines can execute a slowed job.
Beyond that power cap, the algorithm will likely have to kill jobs to fit, regardless of the slowness.
For our machine, those power caps are respectively 0.27 and 0.52.
This explains in part the large number of EcoJobs that were killed in the 0.5 and 0.6 power caps.

\section{Conclusion}
\label{sec:conclusion}

Operating exascale supercomputers requires substantial amounts of power,
which puts a heavy strain in the electricity grid. This strain increases the
risk of power outages, especially in combination with the increased electricity demand
induced by climate events. Facing this risk, exascale supercomputer maintainers
must employ measures to control the power consumption of their platforms, while
keeping the supercomputer operational.

In this paper, we explored the impact of including the users in the decision when it comes to handling a power cap during normal operations of a High-Performance Computing (HPC) supercomputer.
We proposed a scheduling algorithm based on FCFS that can slow jobs through the use of DVFS.
In simulation, we gave users the ability to allow their jobs to be slowed when the platform needs to decrease its power consumption.
In exchange, those users have a lesser chance of seeing their jobs being killed to enforce the power cap.

We performed an experimental campaign in simulation to see the effects of different power caps and different proportions of users that agree to slow their jobs.
Those experiments showed us that, compared to a version of FCFS that only kills jobs when the power cap is enforced, our algorithm has similar performance in throughput and stretch, but is significantly better in number of jobs killed and in use of energy.
Our results are based on the fact that users would volunteer to run their jobs in a
lower power state. One may argue that such a situation may never happen in the
real world, since this willingness to volunteer may not exist in real
users. However, the dynamics of user behavior can be more complex. For instance,
if the proportion of EcoJobs is disclosed to the users:
the higher the proportion of EcoJobs is, the higher the
probability of jobs kills for non-EcoJobs, even though the total number of killed
jobs decreases (Figure~\ref{fig:P_Experiments_Kills-by-scheduler}). This can
potentially influence the decision of the non-EcoJob user. 

Additionally, under a scenario with a significant power cap, an EcoJob
might be slowed or even killed, as illustrated in
Figure~\ref{fig:P_Experiments_Lost-energy-by-scheduler}.
The probability of that event may be high enough that a user may decide it is not worth running EcoJobs. For future work we
plan to better model this user behavior dynamics by using Bayesian decision
theory, for instance. 

We see our algorithm as an add-on that could be fitted to other scheduling algorithms in the same way we did to FCFS in this paper.
As such, we also plan to expand our experiments to other classical scheduling algorithms in HPC, like Easy-Backfilling.

\section*{Acknowledgements}
This work was supported by the REGALE (H2020-JTI-EuroHPC-2019-1 agreement n.
956560), and LIGHTAIDGE (HORIZON-MSCA-2022-PF-01 agreement n. 101107953)
european projects. We also thank Francesco Antici for curating and sharing the
Marconi100 dataset.

\section*{Author contributions}

The authors are listed in alphabetic order.
All the authors participated to the discussions and elaboration of this work. 
Danilo contributed to the data processing,
Luc implemented the methods in Batsim simulation, conducted the experimental protocol and provided experimental results.
Pierre-François helped with implementation and debugging of the methods and the Batsim simulator.
All authors participated into the analysis and results interpretation. 
Danilo was the main writer of Sections~\ref{sec:introduction},~\ref{sec:related-work},~\ref{sec:preliminary-concepts}, and~\ref{sec:methods}.
Luc was the main writer of Sections,~\ref{sec:results}, and~\ref{sec:conclusion}.
Finally, all authors reviewed the final manuscript.

\end{document}